\documentclass{article}
\pdfoutput=1 

\usepackage{graphicx}
\usepackage{amsmath}						
\usepackage{amssymb}
\usepackage{mathtools}
\usepackage{amsthm}
\usepackage{xcolor}

\theoremstyle{definition}

\newcommand{\TheTitle}{Exact Solutions of Time-Delay Integer- and Fractional-Order Advection Equations}

\title{{\TheTitle}}
\providecommand{\keywords}[1]{\textbf{Keywords:} #1}
\providecommand{\msc}[1]{\textbf{2020 MSC:} #1}

\author{
	Christopher N. Angstmann\thanks{School of Mathematics and Statistics, UNSW Australia Sydney, NSW 2052, Australia} \and Stuart-James M. Burney\footnotemark[1]  \and Daniel S. Han\footnotemark[1] \and Bruce I. Henry\footnotemark[1] \and Zhuang Xu\footnotemark[1]}

\begin{document}
	
\maketitle
		
		\begin{abstract}
        Transport phenomena play a vital role in various fields of science and engineering.	
        In this work, exact solutions are derived for advection equations with integer- and fractional-order time derivatives and a constant time-delay in the spatial derivative.        
        Solutions are obtained, for arbitrary separable initial conditions, by incorporating recently introduced delay functions in a separation of variables approach.
        Examples are provided showing oscillatory and translatory behaviours that are fundamentally different to standard propagating wave solutions.
		\end{abstract}

		\keywords{advection equation; time-delay; exact solutions; separation of variables; special functions}
		
		\msc{35C10, 35F10, 34K06, 42A38, 33E20}
		
	\section{Introduction}
Partial differential equations (PDEs) with time-delays arise naturally in modelling multivariable systems where the future state of the system depends not only on the present state but also on the state at an earlier time.
The study of such models is an active area of research in applied mathematics \cite{SYS2015} with modelling applications in  ecology \cite{ECO2008}, engineering \cite{ENG2010}, infectious diseases \cite{GUCOVID} and population dynamics \cite{POP2003}.
%
%
%
%
The reader is referred to the recent monograph by Polyanin, Sorokon and Zurov  \cite{PSZ2024} for an extensive review of exact, approximate and numerical methods for time-delay ordinary differential equations (ODEs) and time-delay PDEs.

Most of the progress to date on time-delay PDEs has focussed on well-known PDEs with a constant delay or a distributed delay in a state variable through an additive reaction term.
Particular examples of this are time-delay advection equations used to model structured cell populations \cite{SYS2015}.
In this work we consider time-delay advection equations where the time-delay appears in the spatial derivative rather than in a reaction term.
These equations can be motivated by fluid dynamics considerations provided that solutions  remain bounded and positive.
However, our interest in time-delay advection equations considered here is not the physical problem of formulating a model equation; it is the applied mathematics problem of finding solutions to a certain class of time-delay PDEs.

In recent work \cite{DF2023}, we provided a systematic introduction to a class of delay functions that can be employed in solving linear integer-order and linear fractional-order delay differential equations. 
Here we have used the method of separation of variables, together with Fourier transform methods and results on delay functions, to provide exact solutions to time-delay advection equations.
The dynamics of the solutions are explored, for particular cases, through examples.
The method of separation of variables has previously been used to compute exact solutions of the time-delay heat equation \cite{PMJ2007,JDM2008}. 
In our work, we have restricted the time-delay to be a constant positive value.
%

%
%
%
%
%

\section{Formulation, Solution and Analysis}	
The time-delay integer-order advection equation that we study here is defined by
\begin{equation} \label{adveceq}
	\frac{\partial u(x,t)}{\partial t}=-a\frac{\partial u(x,t-\tau)}{\partial x}, \quad (x,t)\in (-\infty,\infty)\times (\tau,\infty),
\end{equation}
where $a, \tau \in\mathbb{R^+}$.
This equation can be obtained from the conservation of mass equation
if the advection flux
is time-delayed, but this is not a physical model. 
A more sophisticated form of  Eq.(\ref{adveceq}), which incorporates a delayed Fickian flux and a distribution of delay times, has been proposed as a general model for non-Fickian transport in heterogenous media \cite{DT2006}.
%
%

We first obtain a fundamental solution of Eq.(\ref{adveceq}) subject to the initial condition $u(x,t)=f(x)$ for $(x,t)\in(-\infty,\infty)\times [0,\tau]$ where $f:\mathbb{R}\to\mathbb{C}$ is a smooth and absolutely integrable function.
%
%
This solution will later be used to construct a more general solution subject to a separable initial condition.
We also note that we require $u(x,t)=0$ for $t<0$.
The standard advection equation, with $\tau=0$, can readily be solved with the non-separable initial condition using the method of characteristics.
This leads to the well-known exact propagating wave solution $u(x,t)=f(x-at)$, which signifies a constant horizontal translation of the initial profile of the solution by $at$.
Here we consider $\tau>0$, and we proceed via the method of separation of variables. 
This method can also be used when $\tau=0$, which leads to products of separate trigonometric functions in the variables $x$ and $t$. 
Travelling wave solutions can then be obtained using standard trigonometric identities.

We suppose the solution to Eq.(\ref{adveceq}) is $u(x,t)=X(x)T(t)$, with $X(x)T(t-\tau)\neq 0$, then 
\begin{equation} \label{sepvar}
	\frac{T'(t)}{T(t-\tau)}=-a\frac{X'(x)}{X(x)}=-\lambda, \quad \lambda\in\mathbb{C},
\end{equation}
where $\lambda$ denotes the separation constant. The solution of  the ODE for  $X(x)$ is	\begin{equation} \label{xsol}
	X(x)=c_1(\lambda)e^{\lambda x/a}.
\end{equation} 
It is possible to construct a solution of the ODE for $T(t)$ using Laplace transform methods \cite{PMJ2007,JDM2008}.
Here, we consider an alternative approach using delay functions \cite{DF2023}.
The delay exponential function is defined by the power series 
\begin{equation} \label{def_dexp}
	\mathrm{dexp}(t;\tau)=\sum_{n=0}^\infty\frac{(t-n\tau)^n}{\Gamma(n+1)}\Theta\left(\frac{t}{\tau}-n\right), \quad \frac{t}{\tau}\in\mathbb{R},
\end{equation}
where 
\begin{equation}
	\Theta(y)=\begin{cases} 0 & y<0,\\ 1& y\ge 0, \end{cases}
\end{equation}
is the Heaviside function. It immediately follows from Lemma 1 of \cite{DF2023},
\begin{equation} 
	\frac{d}{dt}\mathrm{dexp}(\lambda t;\lambda\tau)=\lambda\mathrm{dexp}(\lambda (t-\tau);\lambda\tau), \quad  \frac{t}{\tau}\in\mathbb{R}\backslash\{0\},
\end{equation}
that 
\begin{equation} \label{tsol}
	T(t)= c_2(\lambda)\mathrm{dexp}(-\lambda t;-\lambda\tau)
\end{equation}
is a solution of the ODE for the separated solution $T(t)$ in Eq.(\ref{sepvar}), with $T(t)=0$ for $t<0$ and $T(t)=1$ for $0\le t\le\tau$.
Thus, by combining Eqs.(\ref{xsol}) and (\ref{tsol}), we find that a solution to Eq.(\ref{adveceq}) is
\begin{equation}
	u(x,t)=c(\lambda)e^{\lambda x/a}\mathrm{dexp}(-\lambda t;-\lambda\tau).  
\end{equation}
To avoid blow-up, we take $\operatorname{Re}(\lambda)=0$ and without loss of generality we write $\lambda=i2\pi ka$ for $k\in\mathbb{R}$.
The time-delay integer-order advection equation, Eq.(\ref{adveceq}), is linear and homogeneous, so that the principle of superposition applies and we can write the general solution
\begin{equation} \label{supsol}
	u(x,t)=\int_{-\infty}^{\infty} c(k)e^{i2\pi kx}\mathrm{dexp}(-i2\pi ka t;-i2\pi ka\tau)\,dk.
\end{equation}
It remains to find $c(k)$ consistent with the initial condition $u(x,t)=0$ for  $t<0$ and
$u(x,t)=f(x)$ for $0\le t\le\tau$. This can be obtained by employing the initial condition in the above equation, which yields
\begin{equation}
	f(x)=\int_{-\infty}^{\infty} c(k)e^{i2\pi kx}\,dk,\label{fx1}
\end{equation}
since $\mathrm{dexp}(-i2\pi ka t;-i2\pi ka\tau)=1$ for $0\le t\le\tau$. This 
defines $f(x)$ as the inverse Fourier transform of $c(k)$, i.e.
$f(x)=\mathcal{F}^{-1}_k[c(k)](x)$.
By taking the Fourier transform of the above equation and appealing to the Fourier inversion theorem \cite{F1992}, we find that
$c(k)=\hat{f}(k)$. Substituting this quantity into Eq.(\ref{supsol}) gives
\begin{align} 
	u(x,t)&=\int_{-\infty}^{\infty}\hat{f}(k)e^{i2\pi kx}\mathrm{dexp}(-i2\pi ka t;-i2\pi ka\tau)\,dk\\
	&=\int_{-\infty}^{\infty} \hat{f}(k)e^{i2\pi k x}\left(\sum_{n=0}^\infty\frac{(-i2\pi ka)^n(t-n\tau)^n}{\Gamma(n+1)}\Theta\left(\frac{t}{\tau}-n\right)\right)\,dk.
\end{align}
We can interchange the order of integration and summation to now obtain the solution
\begin{align} 
	u(x,t)&=\sum_{n=0}^\infty\left(\int_{-\infty}^{\infty} \hat{f}(k)(i2\pi k)^ne^{i2\pi k x}\,dk\right) \frac{(-a)^n(t-n\tau)^n}{\Gamma(n+1)}\Theta\left(\frac{t}{\tau}-n\right) \nonumber\\
	&=\sum_{n=0}^\infty\frac{d^nf(x)}{dx^n} \frac{(-a)^n(t-n\tau)^n}{\Gamma(n+1)}\Theta\left(\frac{t}{\tau}-n\right)  \label{simplersol}\\
	&=\sum_{n=0}^{\lfloor{\frac{t}{\tau}\rfloor}}\frac{d^nf(x)}{dx^n} \frac{(-a)^n(t-n\tau)^n}{\Gamma(n+1)}. \label{truncsol}
\end{align}
It is simple to verify via direct substitution that the above delay series is indeed a solution of Eq.(\ref{adveceq}) with the initial condition $u(x,t)=f(x)$ for $(x,t)\in(-\infty,\infty)\times [0,\tau]$.

If the $n^{\text{th}}$ derivative of $f(x)$ is independent of $n$ then the factor $d^nf(x)/dx^n$ can be taken outside the sum yielding a solution that has this term as a multiplicative factor of the delay exponential function.
An obvious choice of initial condition in which this occurs is when $f(x)=\exp(x)$.
Note that this is not a physically realistic state on an infinite domain because it is unbounded.
Furthermore, in this special case, it is easy to show that the solution can be negative for sufficiently large $\tau$, which is also an unphysical state.
This solution cannot be used in fluid transport models where $u(x,t)$ represents the concentration of a solute at position $x$ and time $t$.
However, the solution could be useful as a test case for numerical methods.
Additionally, if $f(x)=\exp(x)$, it is known that the solution will oscillate if $\tau>1/ae$ (see \cite{JAD2014} and references therein).
However, this behaviour is not an intrinsic feature of the solution. 
For example, if $f(x)=\exp(-x)$ then it is straightforward to show that the solution is monotonically increasing in time for all $a>0$.

Note that if we consider the limit $\tau\to 0^+$ then Eq.(\ref{simplersol}) becomes
\begin{equation}
	u(x,t)=\sum_{n=0}^\infty \frac{d^nf(x)}{dx^n}\frac{(-at)^n}{\Gamma(n+1)},
\end{equation}
which we recognise as a Taylor series of the function $f(x-at)$ about the point $x$.
Hence,  $f(x-at)$ provides a useful approximation to the solution for $\tau\ll 1$.
The solution reduces to $u(x,t)=f(x)-a(t-\tau)\,df(x)/dx$ for $\tau\le t$ since only the first two terms of the series are nonzero.  

By exploiting the linearity of Eq.(\ref{adveceq}), we can impose an arbitrary separable initial condition to find a broader set of solutions. For instance, a solution of the time-delay integer-order advection equation
\begin{equation} \label{v_adveceq}
	\frac{\partial \tilde{u}(x,t)}{\partial t}=-a\frac{\partial  \tilde{u}(x,t-\tau)}{\partial x}, \quad (x,t)\in (-\infty,\infty)\times (\tau,\infty),
\end{equation}
subject to the initial condition $ \tilde{u}(x,t)=f(x)g(t)$ for $(x,t)\in(-\infty,\infty)\times [0,\tau]$ takes the form 
\begin{equation} \label{ugeneralsol}
	\tilde{u}(x,t)=g(0)u(x,t)+\int_0^\tau g'(s)u(x,t-s)\,ds.
\end{equation}
Recall $u(x,t)$ is the solution given by Eq.(\ref{simplersol}) and $g:\mathbb{R}\to\mathbb{C}$ is a smooth and absolutely integrable function. Again, it is simple to verify via direct substitution that the above equation is indeed a solution of Eq.(\ref{v_adveceq}) with the aforementioned initial condition.

\section{Examples}
%

\subsection*{Example I}
\begin{figure}[h!]
	\centering
	\begin{minipage}{0.493\textwidth}
		\centering
		\includegraphics[width=1\linewidth]{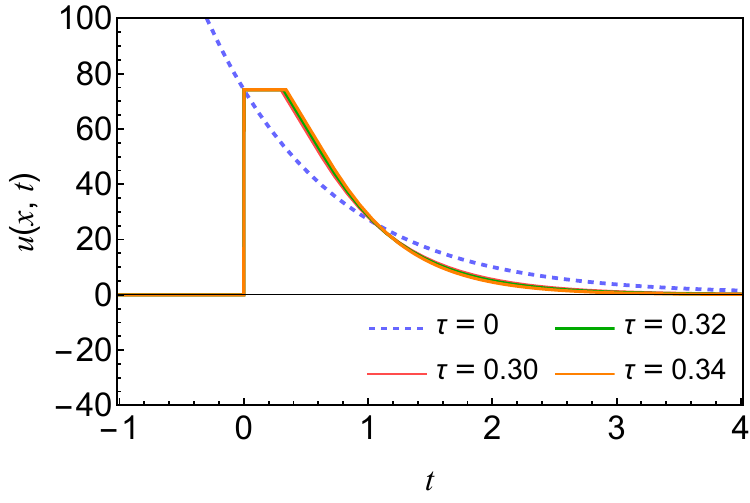}
	\end{minipage}
	\centering
	\hfill
	\begin{minipage}{0.4982\textwidth}
		\centering
		\includegraphics[width=1\linewidth]{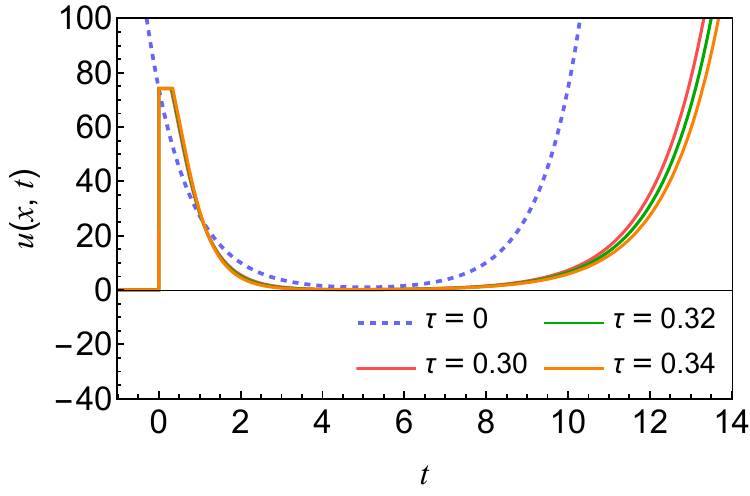}
	\end{minipage}
	\caption{Plot of the solution, Eq.(\ref{sinhsol}), for $\tau=0, 0.30, 0.32$ and $0.34$ with $a=1$ and $x=5$ where $t\in [-1,4]$ (left) and $t\in [-1,14]$ (right). The dashed blue curve is for $\tau=0$. The solid red, green and orange curves are for $\tau=0.30, 0.32$ and $0.34$ respectively.}
	\label{sinhplot1}
\end{figure}

A solution of the time-delay integer-order advection equation with the initial condition $u(x,t)=\cosh(x)$ for $(x,t)\in (-\infty,\infty)\times [0,\tau]$ and $a\in\mathbb{R}^+$ is given by
\begin{equation} \label{sinhsol}
	u(x,t)=\frac{e^x}{2}\mathrm{dexp}(-at;-a\tau)+\frac{e^{-x}}{2}\mathrm{dexp}(at;a\tau).
\end{equation}
The behaviour of this solution is dominated by the first term on the right-hand side as $x$ increases.
For sufficiently large $x$, the solution will exhibit oscillations if $\tau>1/ae$.
The dynamics of the above solution for $\tau=0.30, 0.32$ and $0.34$ are compared with the solution of the standard advection equation, $u(x,t)=\cosh(x-at)$, in Figure \ref{sinhplot1}.
As expected, the above solution is nonoscillatory for these delays when $a=1$.
Instead, it closely approximates the solution of the standard advection equation since $\mathrm{dexp}(-at;-a\tau)\approx e^{-at}$ and $\mathrm{dexp}(at;a\tau)\approx e^{at}$ for $\tau\ll 1$.

\smallskip

\begin{figure}[h!]
	\centering
	\begin{minipage}{0.494\textwidth}
		\centering
		\includegraphics[width=1\linewidth]{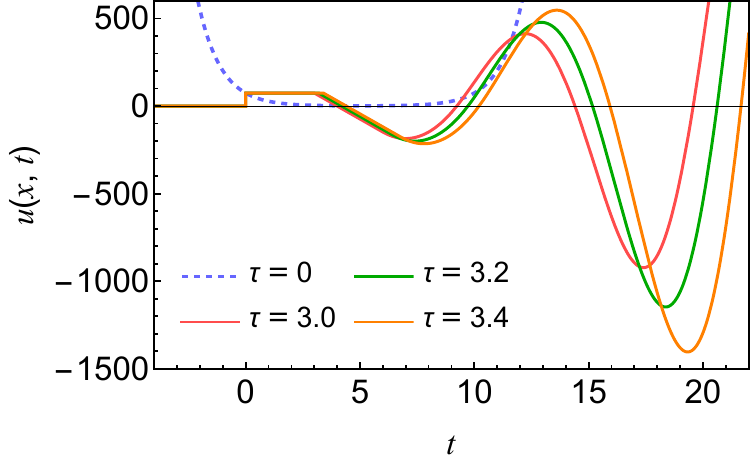}
	\end{minipage}
	\centering
	\hfill
	\begin{minipage}{0.497\textwidth}
		\centering
		\includegraphics[width=1\linewidth]{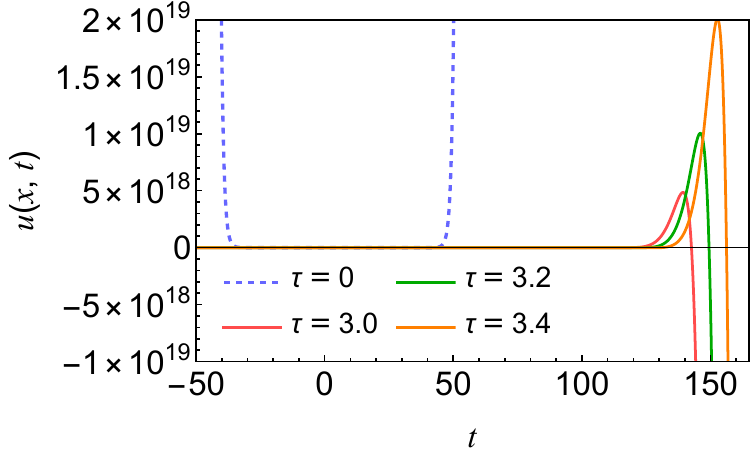}
	\end{minipage}
	\caption{Plot of the solution, Eq.(\ref{sinhsol}), for $\tau=0, 3.0, 3.2$ and $3.4$ where $a=1$ and $x=5$ with $t\in [-4,22]$ (left) and $t\in [-50,165]$ (right). The dashed blue curve is for $\tau=0$. The solid red, green and orange curves are for $\tau=3.0, 3.2$ and $3.4$ respectively.}
	\label{sinhplot2}
\end{figure}

The appearance of oscillations is clear when the delay is increased by an order of magnitude, which is shown in Figure \ref{sinhplot2} for $\tau=3.0, 3.2$ and $3.4$. 
The onset of oscillations occurs more rapidly for larger delays.
Also, the amplitude of the oscillations increases as the delay increases. 
\medskip

\newpage

\subsection*{Example II}
\begin{figure}[h!]
	\centering
	\includegraphics[width=1\linewidth]{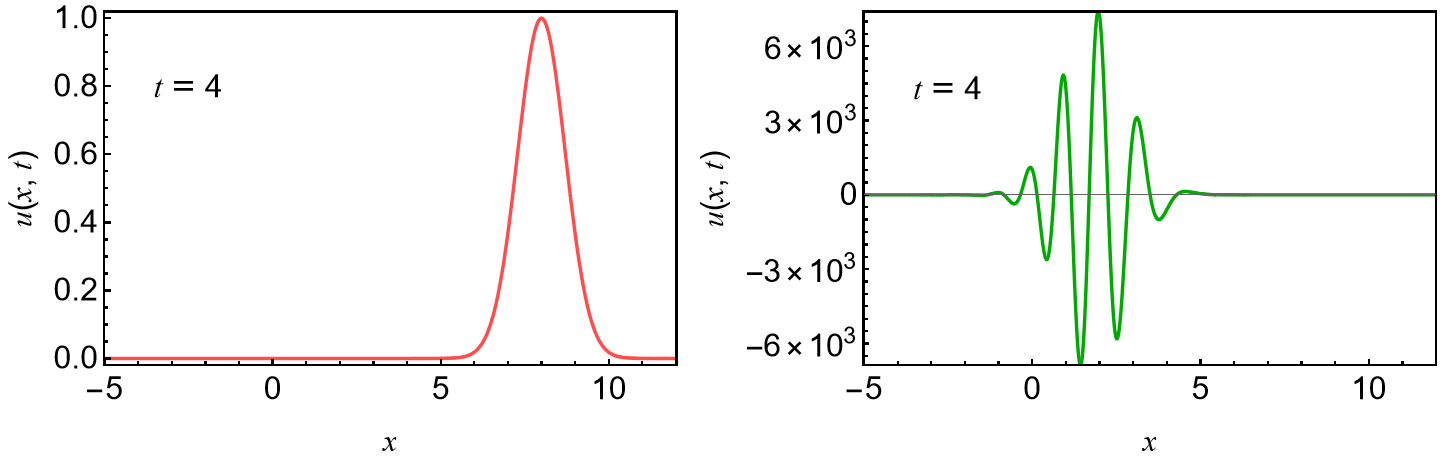}
	\caption{Plots of the solution to the standard advection problem (left) and the time-delay integer-order advection problem (right) when $t=4$ with $a=2$ and $\tau=0.1$.}
	\label{gauss}
\end{figure}

A solution of the time-delay integer-order advection equation with the Gaussian initial condition $u(x,t)=\exp(-x^2)$ for $(x,t)\in (-\infty,\infty)\times [0,\tau]$ and $a\in\mathbb{R}^+$ is given by
\begin{equation} \label{expsol}
	u(x,t)=e^{-x^2}\sum_{n=0}^\infty H_n(x)\frac{a^n(t-n\tau)^n}{\Gamma(n+1)}\Theta\left(\frac{t}{\tau}-n\right),
\end{equation}
where $H_n(x)$ is the $n^{\text{th}}$-order Hermite polynomial.
The solution of the standard advection equation, $u(x,t)=\exp\left(-(x-at)^2\right)$, can be recovered from the above equation in the limit $\tau\to 0^+$ together with the generating function for the Hermite polynomials, 
\begin{equation} 
	\sum_{n=0}^{\infty}H_n(x)\frac{(at)^n}{\Gamma(n+1)}=e^{2axt-(at)^2}. 
\end{equation} 
The dynamics of the above solution, Eq.(\ref{expsol}), are compared with the solution of the standard advection equation in Figure \ref{gauss} when $t=4$.
The Gaussian travelling wave pulse solution of the standard advection equation is replaced by a Gaussian wave packet with increasing amplitude in the time-delay integer-order advection equation.
This could be problematic for using such equations to model physical processes, such as in biological gene regulation \cite{Monk}.

\section{Time-Delay Fractional-Order Advection Equation}
We can generalise the previous results to derive exact solutions to a fractional-order advection equation with time-delay. 
If we replace the classical time derivative in Eq.(\ref{adveceq}) with the Caputo time-fractional derivative of order $\alpha\in (0,1)$ defined by \cite{Caputo1967}
\begin{equation}
	^C\mathcal{D}_t^\alpha v(x,t)\coloneqq\frac{1}{\Gamma(1-\alpha)}\int_0^t \frac{\partial v(x,s)}{\partial s}\frac{1}{(t-s)^\alpha}\,ds,
\end{equation}
then we arrive at the time-delay fractional-order advection equation
\begin{equation} \label{frac_adveceq}
	^C\mathcal{D}_t^\alpha v(x,t)=-a\frac{\partial v(x,t-\tau)}{\partial x}, \quad (x,t)\in (-\infty,\infty)\times (\tau,\infty).
\end{equation}
Note that Eq.(\ref{frac_adveceq}) reduces to Eq.(\ref{adveceq}) by considering the limit $\alpha\to 1$ \cite{Pod}. Similar to above, we consider Eq.(\ref{frac_adveceq}) subject to the initial condition $v(x,t)=f(x)$ for  $(x,t)\in(-\infty,\infty)\times [0,\tau]$. 
The method presented above can be applied here too, noting that the solution of 
$
^C\mathcal{D}_t^\alpha T(t) = - \lambda T(t-\tau)
$
with $T(0)=1$ and $T(t)=0$ for $t<0$, is given by
$
T(t)=\mathrm{dE}_\alpha^{-}\left(-\lambda^{\frac{1}{\alpha}} t; -\lambda^{\frac{1}{\alpha}} \tau\right)
$
where 
\begin{equation}
	\mathrm{dE}_\alpha^-(-t;-\tau)=\sum_{n=0}^\infty (-1)^n\frac{(t-n\tau)^{\alpha n}}{\Gamma(\alpha n+1)}
	\Theta\left(\frac{t}{\tau}-n\right),\quad\frac{t}{\tau}\in \mathbb{R},
\end{equation}
is the delay fractional Mittag-Leffler function \cite{DF2023}.
The solution of Eq.(\ref{frac_adveceq}), with  $v(x,t)=0$ for $t<0$ and $v(x,t)=f(x)$ for $0\le t\le\tau$, is thus given by
\begin{equation} \label{fracsol}
	v(x,t)=\sum_{n=0}^\infty\frac{d^nf(x)}{dx^n} \frac{(-a)^n(t-n\tau)^{\alpha n}}{\Gamma(\alpha n+1)}\Theta\left(\frac{t}{\tau}-n\right).
\end{equation}
Again, we can impose more general forms of the initial condition that permit the solution to vary as a function of time in the interval $[0,\tau]$. 
For instance, a solution of the time-delay fractional-order advection equation
\begin{equation} 
	^C\mathcal{D}_t^\alpha \tilde{v}(x,t)=-a\frac{\partial \tilde{v}(x,t-\tau)}{\partial x}, \quad (x,t)\in (-\infty,\infty)\times (\tau,\infty),
\end{equation}
subject to the initial condition $\tilde{v}(x,t)=f(x)g(t)$ for $(x,t)\in(-\infty,\infty)\times [0,\tau]$ takes the form 
\begin{equation} \label{vgeneralsol}
	\tilde{v}(x,t)=g(0)v(x,t)+\int_0^\tau g'(s)v(x,t-s)\,ds
\end{equation}
where $v(x,t)$ is the solution given by Eq.(\ref{fracsol}).

\section{Summary}
This work serves three purposes.
It shows how delay functions can be incorporated with standard separation of variables to provide solutions of certain time-delay PDEs. 
It shows how time-delays can affect the integrity of solutions where physical requirements need to be met.
It provides exact solutions that can be used for comparisons with other 
methods of solution: exact, approximate or numerical.    
\section*{Acknowledgement}
This research was funded by Australian Research Council grant number DP200100345. We give a special thank you to productive conversations with Boris Huang.
	{
		
	}
\end{document}